\documentclass[11pt,letterpaper]{article}
\usepackage[margin=50pt]{geometry}
\usepackage{amsmath}
\usepackage{physics}
\usepackage{amssymb}
\usepackage{amsfonts}

\usepackage{hyperref}
\usepackage{cleveref}

\usepackage{geometry}
\usepackage{graphicx}
\usepackage[]{subcaption}
\usepackage[inline]{enumitem}

\newtheorem{remark}{Remark}

\newcommand{\tb}{\mathbf{t}}

\def\ccm{Center for Computational Mathematics, Flatiron Institute, Simons Foundation,
  New York, New York 10010}

\def\nyu{Courant Institute of Mathematical Sciences,
  New York University, New York, New York 10012}

\def\kth{Division of Numerical Analysis, Optimization and Systems Theory,
KTH Royal Institute of Technology, SE-100 44 Stockholm}

\def\yale{Applied Mathematics Program,
Yale University, New Haven, CT 06511}
  
\def\mcneel{Robert McNeel \& Associates,
Seattle, WA 98103}
  
\author{
Haiyang Wang\thanks{\nyu\, \& \yale\, 
({\tt hw1927@nyu.edu}).} 
\and
Fredrik Fryklund\thanks{\ccm\, \& \kth\, 
({\tt ffry@kth.se}).} 
\and
Samuel Potter\thanks{\mcneel\,
({\tt sam@mcneel.com}).} 
\and
 Leslie Greengard\thanks{\ccm\, \& \nyu\,
 ({\tt lgreengard@flatironinstitute.org}).}
}
\date{\today}

\title{Scattering theory for Stokes flow in complex branched structures}

\begin{document}

\maketitle

% short version of the abstract
\begin{abstract}
Slow, viscous flow in branched structures arises in many biological and engineering settings. Direct numerical simulation of flow in such complicated multi-scale geometry, however, is a computationally intensive task. 
We propose a ``scattering theory'' framework that dramatically reduces this cost by decomposing networks into components connected by short straight channels. Exploiting the phenomenon of \textit{rapid return to Poiseuille flow} (Saint-Venant’s principle in the context of elasticity), we compute a high-order accurate scattering matrix for each component via boundary integral equations. These precomputed components can then be assembled into arbitrary branched structures, and the precomputed local solutions on each component can be assembled into an accurate global solution.  
The method is modular, has negligible cost, and appears to be the first 
full-fidelity solver that 
makes use of the return to Poiseuille flow phenomenon.
In our (two-dimensional) examples, it matches the accuracy of full-domain solvers 
while requiring only a fraction of the computational effort.
\end{abstract}

% long version of the abstract
% \begin{abstract}
% We introduce a ``scattering theory'' framework that dramatically reduces this cost. The method applies to fairly general geometries composed of components connected by straight channels. The essential assumption is that each connecting channel has length on the order of a few diameters, allowing us to exploit the phenomenon of \textit{rapid return to Poiseuille flow} (Saint-Venant’s principle in the context of elasticity).  
% For each component, we employ a high-order accurate boundary integral solver to construct its scattering matrix, which relates inflow and outflow fluxes to the corresponding pressure drops (in analogy with electromagnetic waveguides). Once precomputed, each component can be reused and combined with others at negligible cost, making the framework highly modular. The flow in a complex network can be assembled at negligible cost using graph theory and the solution of a small sparse linear system. The system size grows linearly with the number of components, much like circuit analysis with Kirchhoff’s laws. To our knowledge, this is the first solver that explicitly leverages the rapid return to Poiseuille flow phenomenon. 
% We present the method in two dimensions and illustrate its accuracy and efficiency with several examples, showing that it achieves results comparable to full-domain solvers at a fraction of the computational expense.
% \end{abstract}

\section{Introduction}

Slow, viscous flow in branched networks arises 
in a number of application areas, 
including hydraulics, microfluidics 
and biological fluid dynamics
\cite{stokes_network_2018,Kirby_2010,langlois,stokes2dexpansion_2018,Pozrikidis_1992,stone2004,stokes2dbifurcation_2025}. 
Modeling such flows accurately is essential both 
for understanding fundamental transport processes and 
for designing engineered devices. 

\begin{figure}
  \centering
  \includegraphics[width=0.8\textwidth]{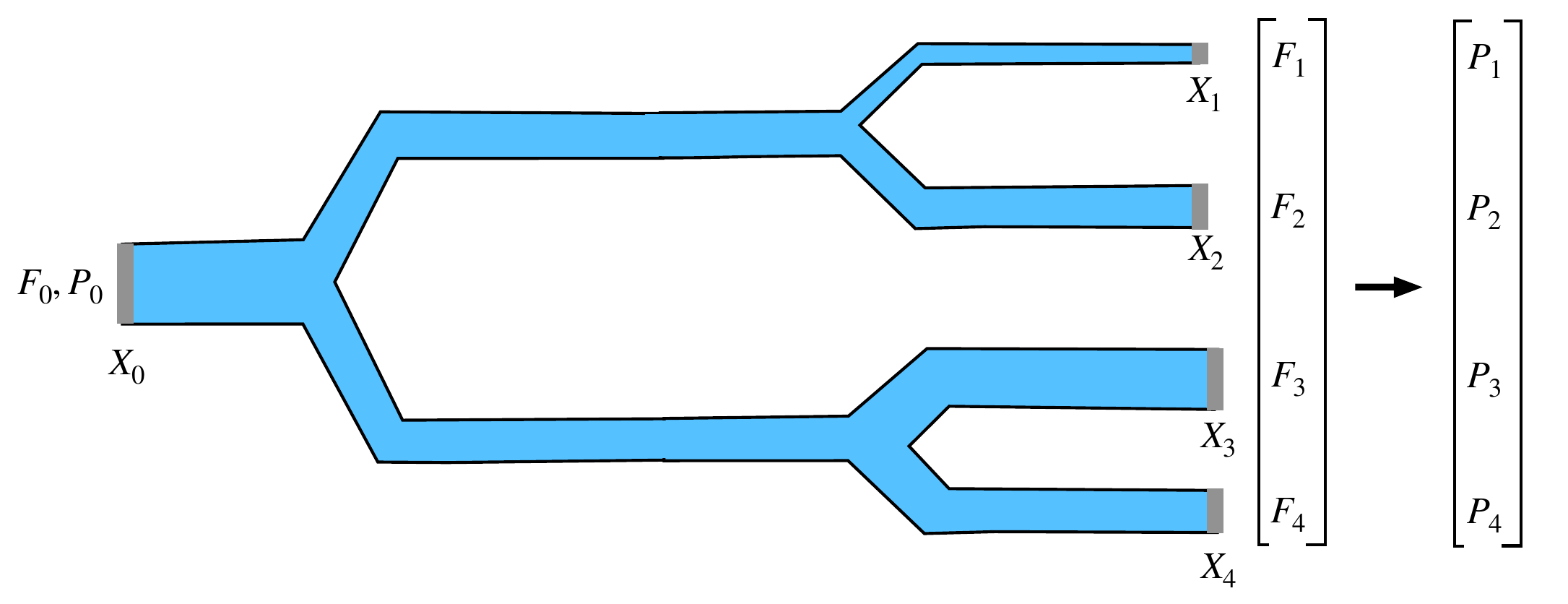}
  \caption{A simple branched network, with inflow at the port $X_0$
and outflows at the ports $X_1,\dots,X_4$.
With specified fluxes $F_1,\dots,F_4$, the flux at the inflow is not an 
independent variable. It must satisfy
$F_0 = F_1 + F_2 + F_3 + F_4$ because the fluid is incompressible.
A straightforward practical question is to determine the necessary 
pressure drops ($P_1-P_0,\dots,P_4-P_0$) 
from $X_0$ to each of the outflow ports in order to establish
the desired flow profile. ($P_0$ can be set arbitrarily since the pressure
is only defined up to an arbitrary constant.)
This requires solving the Stokes equations with the prescribed velocities.
A more complex optimization question would be to \emph{design}
the interior structure of a network with specified inflows and outflows in order to
minimize the pressure drops needed to establish the desired flow profile.
This would require solving a new boundary value problem for each proposed geometry 
within an outer design loop.}
\label{design_fig}
\end{figure}

For an incompressible fluid in two dimensions 
in the zero Reynolds number limit, 
the governing equations are those of Stokes flow:
\begin{align}
  \mu \Delta u = \pdv{p}{x},\quad & \mu \Delta v = \pdv{p}{y}
  \label{stokes}                                                                       \\
  \pdv{u}{x} + \pdv{v}{y}                      & = 0
  \label{continuity}
\end{align}
where $u,v$ are components of velocity, $p$ is the pressure,
and $\mu$ is the dynamic viscosity. 
The corresponding vorticity is defined as
$\zeta  = u_y - v_x.$ 
A standard problem for Stokes flow is to specify boundary 
conditions on the velocity and to determine the velocity and pressure at all
other points in the domain of interest
\cite{langlois,Pozrikidis_1992,ladyzhenskaya}: 
We will restrict our attention here to flows in structures of the type depicted
in \cref{design_fig}, a branched network of channels with no-slip boundary conditions
on the channel walls ($u = 0$) and specified fluxes at the various inlets/outlets,
which we will refer to as \emph{ports}.

\begin{figure}
  \centering
  \includegraphics[width=0.8\textwidth]{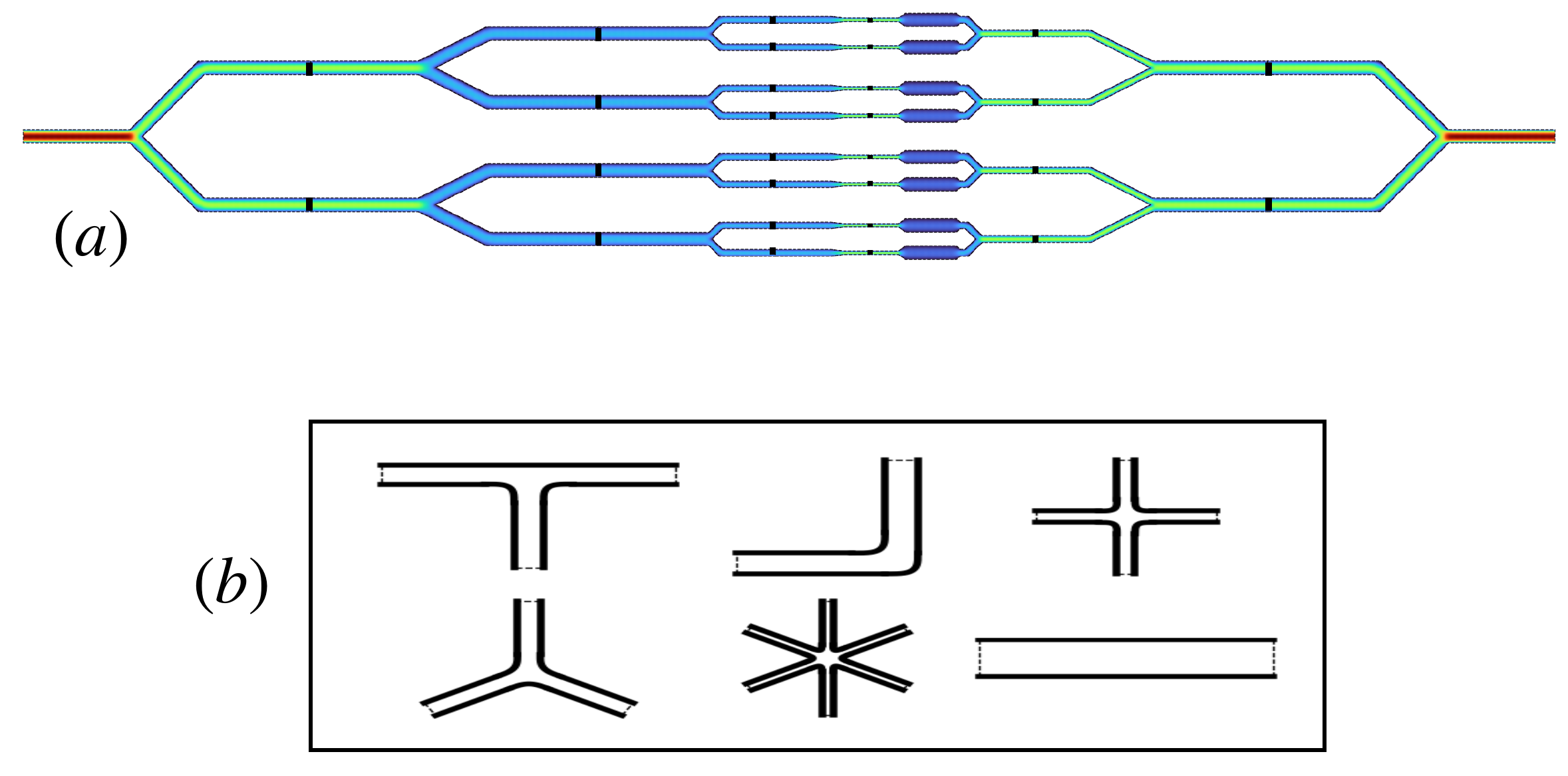}
  \caption{(a) A more elaborate branched network  and (b) a selection of possible
components. These components can be coupled at straight channel interfaces
if they have the same cross-section (diameter) and are
separated from complex features by a few channel widths,
such as the locations marked by vertical black bars in (a).}
\label{fig-components}
\end{figure}

In the last decades, linear-scaling fast algorithms have been developed to solve
\eqref{stokes}, \eqref{continuity} using boundary integral 
formulations, described very briefly in \cref{sec_inteq}.
However, applying such techniques to the full branched networks is 
computationally intensive: each new geometry requires solving a large-scale boundary value problem, often with millions of degrees of freedom and intricate multiscale boundaries. 
This motivates the search for a more modular, reduced-order approach.
Here,
we propose a scattering-matrix based framework that circumvents the cost of full-domain solvers. First, we break the branched network into {modular components} (see \cref{fig-components}). Using high-order integral equation solvers, we compute the \textit{scattering matrix} of each component, which maps fluxes at ports to the associated pressure drops. 
Once computed, these matrices can be reused to assemble large networks at negligible cost, using matching conditions, a little graph theory, and the solution of 
a small sparse linear system. 
Second, and somewhat surprisingly, our method appears to be the first to systematically
exploit the phenomenon of \textit{rapid return to Poiseuille flow} in a 
high-order accuracy manner. 

It is well known that when fluid enters a straight channel, it converges exponentially 
fast to what is called Poiseuille flow within just a 
few diameters~\cite{gregoryTractionBoundaryValue1980,rogerson,horganDECAYESTIMATESBIHARMONIC1989} (see~\cref{sec:ret2poi}), and that
matching Poiseuille profiles at such points is trivial.
To guarantee high-order accuracy, this imposes a mild restriction on the allowed
geometry – namely, that the component library from which the
overall structure is built contains a short span of a straight channel at each port, as 
illustrated in \cref{fig-components}.

We should note that the rapid return to Poiseuille flow \emph{has} been taken advantage of 
before, as has the analogy between Stokes flow in devices with branching
channels and Kirchhoff's laws in electrical networks
(see, for example
\cite{stokes_network_2018,marusic,stephenson2024}).
For this, fluid flux corresponds to current, pressure to voltage and the scattering matrix
to resistors.
However, in order to obtain accurate flows 
in the full device, the solution to the Stokes equations 
must be computed in each component with full fidelity and the matching imposed 
at points where the flow has a simple Poiseuille flow profile. 
In this paper, we use high-order accurate numerical methods to characterize
the scattering matrices for each component and an outer framework for
coupling many such components together.

The result is a modular and efficient solver that achieves the accuracy of 
full-domain Stokes simulations at a fraction of the computational expense. 
The remainder of the paper is organized as follows. In the next section,
we briefly present the mathematical foundations of the method,
and the classical analysis of flow in a semi-infinite
channel (the return to Poiseuille phenomenon).
In \cref{sec:scatmat}, we show how to construct the scattering
matrix for a single component and in \cref{sec:merge},
we show how to solve the global coupling problem.
This includes a brief summary of the graph
theory issues needed to resolve situations when there are loops
present in the network and recirculation must be taken into account.
We present numerical results in \cref{sec:numericalresults} 
and discuss avenues for further investigation in 
\cref{sec:conclusions}.

% More precisely, we use a boundary integral equation method to characterize
% the behavior of each component as a ``scattering matrix," 
% relating the pressure drop
% across inflow and outflow channels with the fluid flux (adopting the
% language from its use in electromagnetic and acoustic applications).

\section{Mathematical Preliminaries\label{mathprelim}}

For our purposes, it is sufficient to
consider the Stokes equations in
a bounded domain $D\subset \mathbb R^2$, 
with boundary $\partial D=\Gamma$.  
We assume that the velocity $\mathbf{u} = (u,v)$ is specified on
the boundary $\Gamma$, with
\begin{equation}
u(\tb) = h_2(\tb),\quad v(\tb) = - h_1(\tb), \quad \tb \in \Gamma.
  \label{bdr-velocity}
\end{equation}
The boundary value problem \eqref{stokes}, \eqref{continuity} 
subject to \eqref{bdr-velocity} is well-known to have a unique 
solution 
\cite{kim2005microhydrodynamics,ladyzhenskaya,mikhlin1964integral,Pozrikidis_1992},
up to a constant value for the pressure. A convenient reformulation of Stokes flow uses the scalar stream function $W(x,y)$ with 
\begin{equation}
  \pdv{W}{x} = -v,\quad \pdv{W}{y} = u,  \quad 
\Delta W = \zeta \label{stream-1}.
\end{equation}
The incompressibility condition \eqref{continuity} is 
satisfied by construction, and it is easy to verify that the stream function satisfies the biharmonic equation, with the boundary condition \eqref{bdr-velocity} translating to \eqref{bih-bv}: 
\begin{align}
& \Delta^2 W(x,y) = 0,   & (x,y)\in D \label{biharmonic} \\
   & \pdv{W}{x}(\tb) = h_1(\tb),\quad \pdv{W}{y}(\tb) = h_2(\tb),  & \tb \in \Gamma. \label{bih-bv}
\end{align}

Thus, solving the Stokes equations with prescribed boundary velocities is equivalent to solving the biharmonic boundary value problem for the stream function $W$. In the next subsection, we will briefly discuss the numerical solution of the biharmonic equation.

\subsection{Numerical Solution of the Biharmonic Equation} \label{sec_inteq}
A variety of numerical methods can be used to 
solve \cref{biharmonic,bih-bv}, 
including conformal mapping~\cite{delillo}, 
rational approximation~\cite{lightning_stokes},
or boundary integral equation formulations~\cite{greengardIntegralEquationMethods1996}. Among these, integral equation methods are particularly effective as they reduce the problem dimension by one. The literature on integral equations for Stokes flow
is substantial, and we refer the reader to the texts
\cite{kim2005microhydrodynamics,
ladyzhenskaya,
mikhlin1964integral,
Pozrikidis_1992} 
as well as the papers
\cite{Davis1985,
GREENBAUM1992216,
greengardIntegralEquationMethods1996,
HelsingJiang,
mbgv,
OJALA2015145,
rachhSolutionStokesEquation2020,
wuSolutionStokesFlow2020} 
and the references therein. 

In this work, we use the method of~\cite{greengardIntegralEquationMethods1996},
following the classical treatment of Goursat, Mikhlin, Muskhelishvili and others
\cite{goursat,mikhlin1964integral,muskhelishviliBasicProblemsMathematical1977},
based on the fact that
any plane biharmonic function $W(x,y)$ can be expressed in the form:
\begin{equation}
  W(x,y) = \Re (\bar z \phi(z) + \chi (z)) \label{Goursat}
\end{equation}
where the Goursat functions $\phi, \chi$ are analytic functions of 
the complex variable $z = x+iy$
and $\Re(f)$ denotes the real part of the function $f$.
With a slight abuse of notation,
we identify $(x,y) \in \mathbb{R}^2$ with $z=x + iy \in \mathbb{C}$.
A simple calculation shows that the
velocity field, vorticity and pressure are given by 
\begin{align}
  \pdv{W}{x} + i\pdv{W}{y}
    &= \phi(z) + z \overline{\phi'(z)} + \overline{\psi(z)}
    \label{muskhelishvili}\\
  \zeta + \frac{i}{\mu}p & 
    = 4\phi'(z) \label{pressure-and-vorticity}
\end{align} 
where $\psi = \chi'$.
\eqref{muskhelishvili} is often referred to as 
the Muskhelishvili formula. Using Goursat functions,
the boundary value problem~\eqref{biharmonic} and~\eqref{bih-bv} can be rewritten as
\begin{align}
  \phi(t) + t\overline{\phi'(t)} + \overline{\psi(t)} = h(t), 
  \quad t \in \Gamma \label{musk-bvp}
\end{align} where $h(t) =  h_1(t) + i h_2(t)$,  
and $t$ is understood as a complex variable.

The classical integral equation for the biharmonic equation
in multiply-connected domains is 
due to Sherman and Lauricella and (for Stokes flow) based on the 
ansatz~\cite{greengardIntegralEquationMethods1996}:
\begin{align}
  \phi(z) & =
  \frac {1}{2\pi i} \int_\Gamma \frac{\omega(\xi)}{\xi - z} d\xi
  + \sum_{k=1}^M C_k \log (z-z_k) \label{sl-phi}
  \\
  \psi(z) & =
  \frac {1}{2\pi i} \int_\Gamma \frac{\overline{\omega(\xi)}d\xi  +  \omega(\xi)\overline{d\xi}}{\xi - z} \\
  & \quad - \frac {1}{2\pi i} 
    \int_\Gamma \frac{\overline{\xi} \omega(\xi)}{{(\xi - z)}^2} d\xi  \label{sl-psi} \nonumber \\
          & \quad + \sum _{k=1}^M
  \left( \frac{b_k}{z-z_k} + \overline C_k \log (z-z_k) -  C_k \frac{\overline z_k}{z-z_k} \right) \nonumber
\end{align}
where $\omega$ is an unknown complex density on $\Gamma$,
$z_k$ are arbitrary prescribed points inside the component curve $\Gamma_k$,
and $C_k, b_k$ are constants defined by:
\begin{align}
  C_k = \int_{\Gamma_k} \omega(\xi) |d\xi|, \quad b_k = 2 \Im\int_{\Gamma_k} \overline{\omega(\xi)} {d\xi}
\end{align}

It is worth noting that, while $\psi$ and $\phi$ are multi-valued, the physical
quantities of interest (velocity, pressure, and vorticity) are all single-valued 
functions of $z$. We refer the reader to~\cite{greengardIntegralEquationMethods1996}
for a discussion of discretization and the use of fast multipole acceleration.
The choice of integral equation solver is not the central contribution of this paper. 
We should, however, note the work of Luca and Llewellyn Smith 
\cite{stokes2dexpansion_2018} which couples
the complex variable formulation 
with the unified transform method of Fokas. The goal of that work is closely related
to ours:
the characterization of the pressure head required to achieve a given flow profile in
a channel with a linear expansion. See also \cite{crowdy2013}.

\subsection{Poiseuille flow \label{sec:ret2poi}}

An important particular solution to the Stokes equations in an infinite channel
$D = \{ (x,y) | x \in \mathbb{R}, \  -L \leq y \leq L \}$ is that of Poiseuille flow:
a constant pressure gradient and a parabolic, laminar velocity profile:
\[
\nabla p = (C,0), u = C \, (y-L)(y+L)/(2\mu), v = 0 .
\]
 The flux $F$ in a channel at a specific location $x=x_0$ is defined as
\begin{equation}
 F = \int_{-L}^L u(x_0,y) dy,
\label{fluxdef}
\end{equation}
and, in the case of Poiseuille flow, is simply
$F = 2CL^3/(3\mu)$.

A remarkable feature of Stokes flow is that 
any velocity profile entering a straight channel  
rapidly converges to Poiseuille flow to many digits of accuracy after only a 
few channel diameters
\cite{gregoryTractionBoundaryValue1980,coRecentDevelopmentsConcerning1983,horganDECAYESTIMATESBIHARMONIC1989,knowlesENERGYESTIMATEBIHARMONIC,stokes2dexpansion_2018,rogerson}.
It is this fact, which we refer to as 
\emph{rapid return to Poiseuille flow},
that permits the use of simple scattering matrices 
with high accuracy.

\begin{figure*}[ht]
  \centering
  \includegraphics[width=0.7\textwidth]{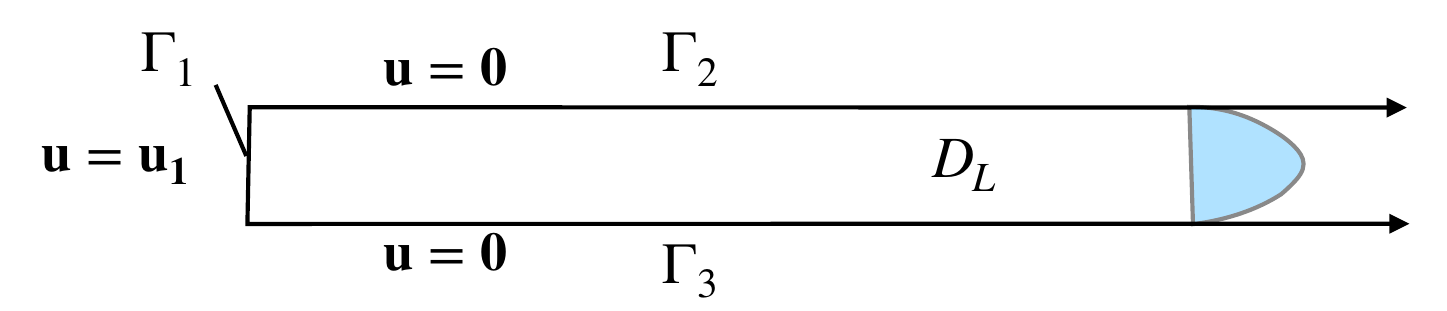}
  \caption{
    A semi-infinite channel in the $x$-direction,
    with arbitrary inflow profile imposed on the left
    ${\bf u} = {\bf u}_1$, and zero velocity conditions on the top and bottom walls.
    \label{fig:return}}
\end{figure*}

The theoretical analysis of this phenomenon can be found in
\cite{gregoryTractionBoundaryValue1980,rogerson}. In brief,
consider a semi-infinite straight channel of width $2L$: 
$D_L = \{(x,y)\mid x \ge 0, -L \le y \le L \}$
with boundary: 
\begin{align}
\Gamma_L & = \Gamma_L^1 \cup \Gamma_L^2 \cup \Gamma_L^3 \\ 
& =\{(0,y)| -L \leq y \le L \} \cup \{(x,L)|\mid x\ge 0\} \cup 
\{(x,-L)\mid x\ge 0\} \, . \nonumber
\end{align}
On the top and bottom walls, $\Gamma_L^2,\Gamma_L^3$, we impose
zero velocity conditions ${\bf u} = (0,0)$.
On the inlet $\Gamma_L^1$ we assume an arbitrary velocity profile
${\bf u}_1 = (f,g)$, assuming only that $f(0,L) = f(0,-L) = g(0,L) = g(0,-L) = 0$
(see \cref{fig:return}).

%\begin{definition} \label{fluxdef}
%The flux at a distance $x$ along the channel is defined as
%$F = \int_{-L}^L {\bf t} \cdot {\bf u}(y) dy$, where ${\bf t}$ is parallel
%to the channel direction.
%\end{definition}

This problem has a semi-analytic solution, 
so long as 
$f^{\prime\prime},g^{\prime\prime}$ exist and are of bounded 
variation~\cite{crowdy2013,gregoryTractionBoundaryValue1980,horganDECAYESTIMATESBIHARMONIC1989,karageorghis89,coRecentDevelopmentsConcerning1983,rogerson,smith}. 
More precisely, the stream function for $x > 0$ takes the form
\begin{equation}
W(x,y) = 
C_0 \, (y^3/3- y \, L^2)/(2\mu)
+ \sum_{n=1}^\infty
C_n \phi_n(y) e^{-p_n x/(2L)},
\label{papkovich}
\end{equation}
where the eigenvalues $p_n$ are the complex roots of $\sin^2(2z) - 4z^2 = 0$ with 
positive real parts. The functions $\phi_n$ in the infinite series 
are known as Papkovich-Fadle eigenfunctions. 
The first term corresponds to Poiseuille flow and
the decay rates of the higher order corrections are governed by the real parts
of the $p_n$. The slowest decaying mode has $\text{Re }p_1 \approx 4.2$, 
implying that non-Poiseuille components decay at a rate of $\exp(-4.2x/L)$. 
Thus, after only $4$ channel diameters, the flow matches Poiseuille to roughly 7 digits, and after $8$ diameters to about 14 digits \cite{gregoryTractionBoundaryValue1980,rogerson,horganDECAYESTIMATESBIHARMONIC1989,smith}.

\begin{figure}[!h]
\centering
\includegraphics[width=0.7\textwidth]{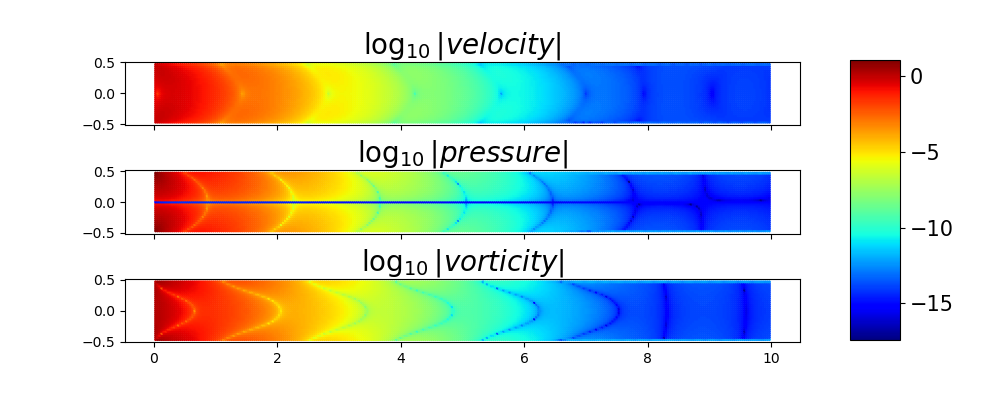}
\caption{The magnitude of the velocity, pressure and vorticity
in a straight channel, with a random inlet velocity that has no Poiseuille
flow component, plotted on a logarithmic scale. \label{fig:rtppipe}}
\end{figure}

To verify this numerically,
we solve the Stokes equations in a straight channel of width $1$ 
and length $10$ that satisfies zero 
velocity boundary conditions on the top, bottom, and right boundaries. 
At the left inlet, we choose a random velocity field with zero flux and
measure the rate at which the flow dissipates to zero in the downstream direction.
In terms of the representation \eqref{papkovich}, this corresponds to a biharmonic
function with no projection on the Poiseuille component of flow, permitting 
us to assess the rate of decay of the Papkovich-Fadle eigenfunction components
themselves.
In \cref{fig:rtppipe} we plot the magnitude of the velocity, pressure and 
vorticity along the channel, showing excellent agreement with the analytic
estimate that the decay rate is $e^{-4.2x}$.
(To check that our numerical solution is converged in the semi-infinite setting,
we repeat the calculation for a channel of length $20$ and find agreement
to machine precision.)

\section{Computing the scattering matrix for a single component} \label{sec:scatmat}

\begin{figure*}[ht]
  \centering
  \includegraphics[width=0.7\textwidth]{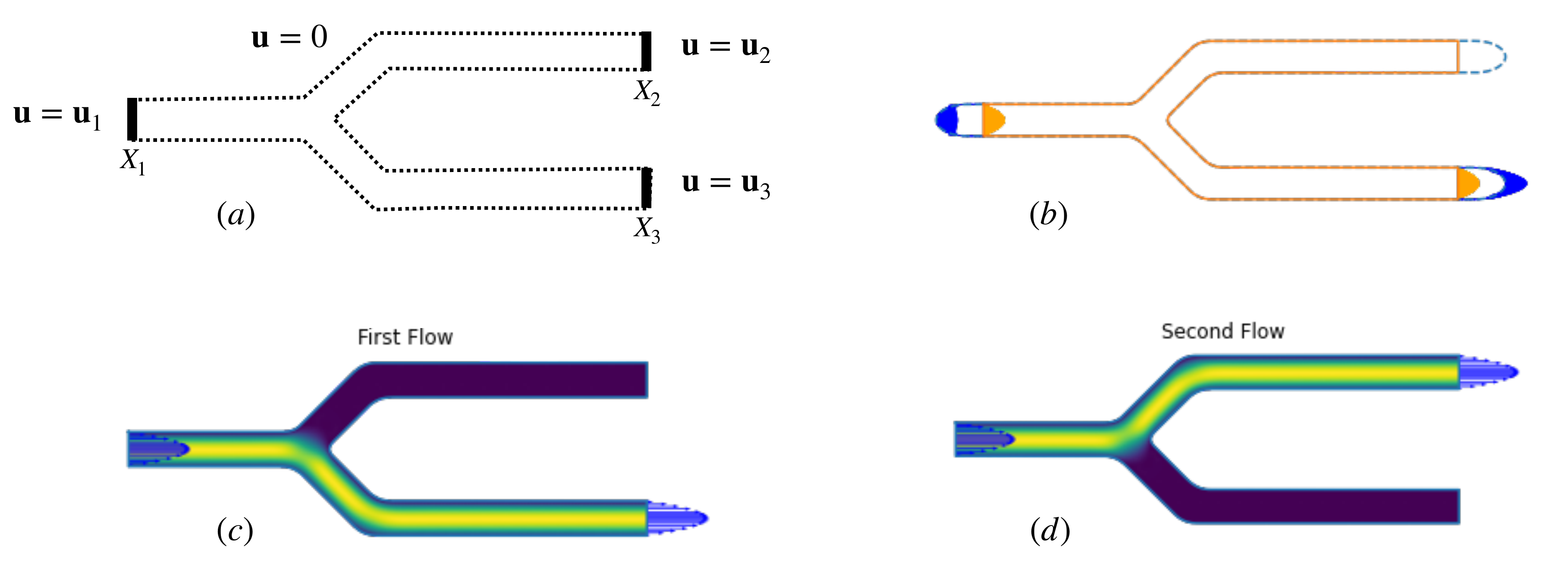}
  \caption{
    (a) A $Y$-junction with zero velocity boundary conditions on the walls (dashed lines)
 and specified velocities ${\bf u}_1, {\bf u}_2, {\bf u}_3$ at the three inlets/outlets.
(b) The interior corners are locally smoothed and resolved with refinement to the scale
of the curvature. Rather than imposing Poiseuille flow on the vertical segments in
(a), it is imposed on smooth caps (dashed blue lines). 
Two boundary value problems are then solved. First, we establish a Poiseuille profile
at the inlet and lower outlet with fluxes $F_1=F_3=1$, and ${\bf u}_2=0$, illustrated in (c).
Second, we switch the conditions on the outlets, with Poiseuille flows
at the inlet and upper outlet with $F_1=F_2=1$ and ${\bf u}_3=0$,
illustrated in (d). }
\label{fig:bvp}
\end{figure*}

Let us first consider
a $Y$-junction as shown in \cref{fig:bvp} (a).
The Stokes equations are solvable so long as the total flux is zero,
and we use the Sherman-Lauricella equation to compute the solution,
from which we obtain the corresponding pressure profile.

While high order accuracy can be maintained, even with corner singularities
~\cite{wuSolutionStokesFlow2020,rachhSolutionStokesEquation2020}, in this paper
we have used only smooth boundary curves for the sake of simplicity. 
Starting with a polygonal domain,
we smooth the interior corners in a small neighborhood of the singularity using
the method of \cite{epsteinSmoothedCornersScattered2016}.
Smooth caps are placed on the inlets and outlets
following the construction of \cite{baggeHighlyAccurateSpecial2021}, illustrated in
\cref{fig:bvp} (b). For this $Y$-junction, we solve
two boundary value problems: 

\begin{description}
\item (A) with Poiseuille velocity profiles at $X_1,X_2,X_3$ and 
fluxes $F_1=1, F_2=0, F_3=1$, 
(\cref{fig:bvp} (c)),  
\item (B) with Poiseuille velocity profiles at $X_1,X_2,X_3$ and 
fluxes $F_1=1, F_2=1, F_3=0$ (\cref{fig:bvp} (d)). 
\end{description}

After solving the Sherman-Lauricella equation, we have the two ``basis'' flows 
(flow $A$ and flow $B$).
Let us denote the computed pressure drops
from inlet $X_1$ to the respective outlets by $p_{2A}, p_{3A}$ 
for flow $A$ and by $p_{2B}, p_{3B}$ for flow $B$.
By the linearity of the Stokes equations, 
it follows that for any $F_2$ and $F_3$ (with $F_1 = F_2+F_3$ to satisfy the zero
net flux condition), the pressure drops from inlet 1 to outlet 2 and from inlet 1 to outlet 3 must be
$p_2 = F_2 p_{2A} + F_3 p_{2B}$ and 
$p_3 = F_2 p_{3A} + F_3 p_{3B}$.
In matrix form, 
\[ 
\begin{pmatrix}
p_2 \\ p_3
\end{pmatrix}  =  S 
\begin{pmatrix}
F_2 \\ F_3
\end{pmatrix}, \qquad
{\rm with}\ 
S =
\begin{pmatrix}
p_{2A} & p_{2B} \\ 
p_{3A} & p_{3B} 
\end{pmatrix} .
\]
We refer to the matrix $S$ as the \emph{scattering matrix} for the component.
The extension of the scattering matrix to a general 
component with a total of $m$ ports (inlets/outlets) is straightforward.
$S$ is of dimension $m-1$, since there are only $m-1$ independent
fluxes if we are to maintain the zero total flux condition.
Fixing the first port as the inlet, labeled ``1'' to enforce the condition 
$F_1 \equiv \sum_{j=2}^{m} F_j$,
we may generate column $j-1$ of the scattering matrix 
by solving the boundary value problem with $F_1 = F_j = 1$ and all other fluxes
set to zero, and computing the $(m-1)$ pressure drops with respect to the first port.

\section{Assembling the global flow} \label{sec:merge}

Given the scattering matrices for the various components, the global solution is
obtained by concatenating the fluxes so that they match at the
component interfaces. When loops are present (that is, the domain is not
simply-connected), we also require that the pressure be a continuous function 
in the domain. (In DC electrical networks, these are analogous to
Kirchhoff's current law and Kirchhoff's voltage law, respectively.)

\begin{figure}[ht]
  \centering
  \includegraphics[width=0.7\textwidth]{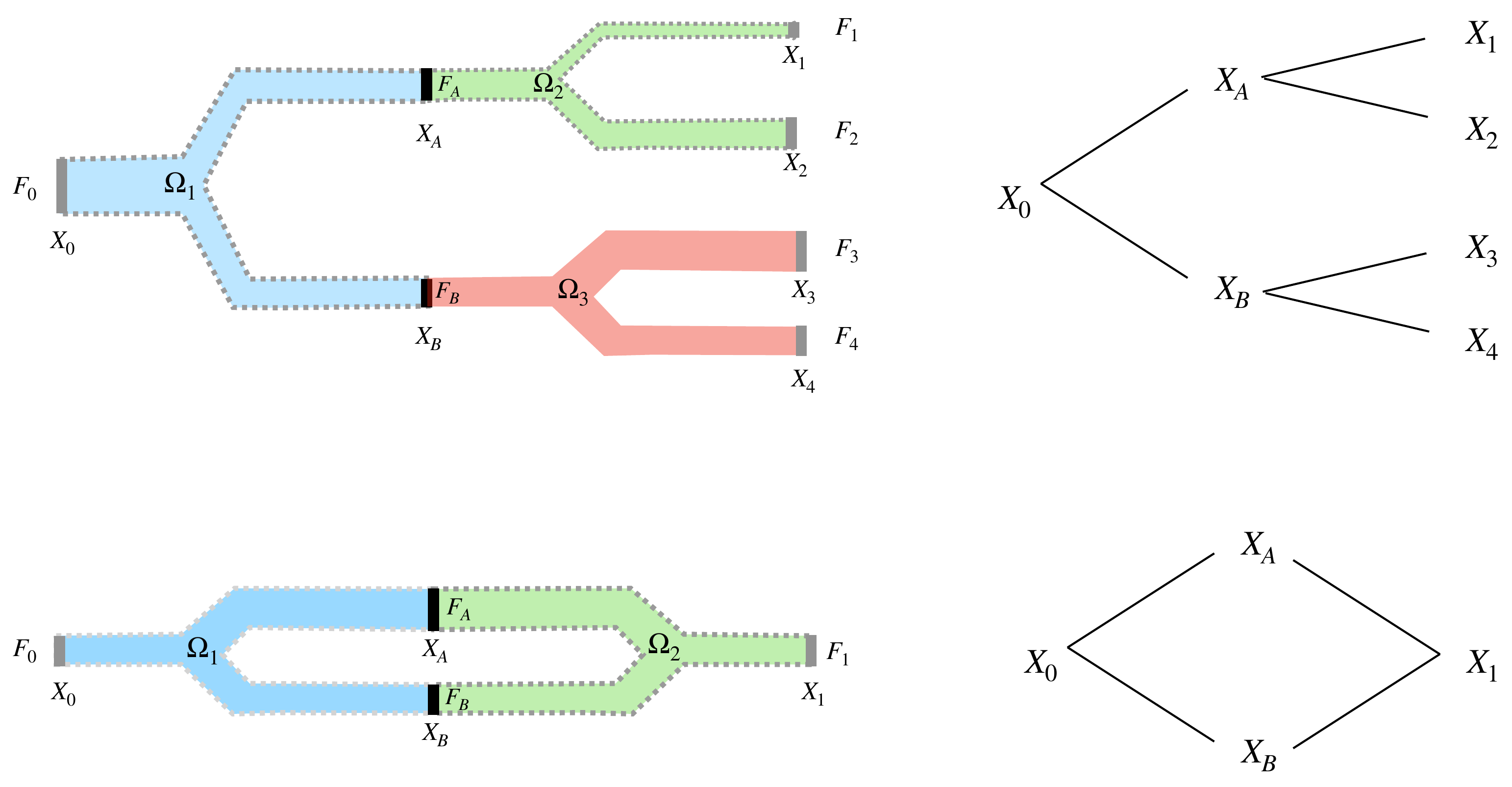}
  \caption{
    In the upper panels we consider coupling three different $Y$-junctions together,
assuming only that, where they join, the channel widths are identical.
The matching condition is applied at the locations marked by thick vertical bars,
where the fluxes $F_A$ and $F_B$ are initially unknown. 
In the graph on the right, $X_0$ corresponds to the location of the inlet,
$X_A$ and $X_B$ to the matching interface and 
$X_1$, $X_2$, $X_3$, $X_4$ to the four outlet locations.
In this diagrammatic abstraction of the problem we need to find the fluxes
$F_A$, $F_B$, given $F_1,\dots,F_4$ with $F_0 = \sum_{j=1}^4 F_j$.
From the zero net flux condition, it is clear that $F_A = F_1+F_2$ and that
$F_B = F_3+F_4$. Once all fluxes are known, the scattering matrices for the
individual components can be used to determine all pressure drops. 
Note that no linear system needs to be solved.
In the lower panels we consider coupling two $Y$-junctions together,
with the result that there is only one input port and one output port.
In this case the domain $(\Omega_1 \cup \Omega_2)$ is not simply connected.
While $F_A+F_B = F_1 = F_0$, we need another linearly independent condition 
to solve for $F_A,F_B$. That condition is that the pressure is continuous.
More precisely $p_{0,B}+p_{B,1}+p_{1,A}+p_{A,0} = 0$, where 
$p_{0,B}$ denotes the pressure drop from $V^0$ to $V^B$, etc.
}
\label{fig:merge}
\end{figure}

To illustrate how these conditions are used in practice,
let us consider the two examples in \cref{fig:merge}.
In the upper panels, three $Y$-junctions are connected, with outflows
$F_1,F_2,F_3,F_4$ specified and $F_0 = F_1+F_2+F_3+F_4$.
The domain in this case is simply connected, so that 
the corresponding graph is acyclic (shown on the right, 
where the interior nodes represent the matching
interfaces). Beginning with the ``leaf nodes'' of the graph, we can recursively
determine the fluxes at all interior nodes. Once the fluxes are known, the 
scattering matrices for the individual components determine all desired
pressure drops.
In the lower panels, two $Y$-junctions are coupled to create a doubly-connected
domain. Since $F_A$ and $F_B$ are both unknown, they cannot be determined 
from the single input flux $F_0=F_1$.
Because of the loop present in the full device geometry (the cycle in the 
corresponding graph), however, we also need to enforce the continuity of pressure when the 
loop is traversed. Letting $p_{\alpha,\beta}$ denote the pressure drop from node
$\alpha$ to node $\beta$, we must have
$p_{0,B}+p_{B,1}+p_{1,A}+p_{A,0} = 0$. Letting  $S_1$ and  $S_2$ denote the
$2 \times 2$ scattering matrices for $\Omega_1$ and $\Omega_2$, 
it is straightforward to see that
\[ 
[S_1(2,1)- S_2(2,1) + S_2(1,1)- S_1(1,1)] F_A +
[S_1(2,2)- S_2(2,2) + S_2(1,2)- S_1(1,2)] F_B = 0.
\]
This equation, together with the fact that $F_A+F_B=F_0$ permits us to solve for
the two unknown fluxes.

The general case involves the imposition of flux continuity at all interior nodes
and the continuity of pressure condition applied to all independent cycles.
We will refer to this as the \emph{assembly matrix}.
It can be handled by a mixture of graph theory and linear algebra.
For an acyclic device, one simply traverses the tree from leaf nodes to the root,
computing fluxes at all interior nodes from the zero net flux condition.
All pressure drops can then be computed from the individual scattering matrices.

When cycles are present, it can be shown that the discrete problem is solvable if the
zero net flux condition is applied at all interior nodes, and continuity
of pressure is added for each independent cycle
(as shown in the simple case above) \cite{feldmann}.
Since this linear algebraic aspect of the problem
is the same as Kirchhoff's current and voltage laws and there is a substantial 
literature discussing circuit theory, we refer the reader to the literature
\cite{clayton,desoer,gotlieb,paton} and provide open source software for this 
purpose \cite{stokesnetwork}.

\begin{remark}
For modest-sized networks with cycles, 
direct Gaussian elimination has negligible cost. 
For large-scale problems, more sophisticated methods can be used, such as
multifrontal sparse direct solvers \cite{davis_direct,duffreid}, the recursive
merging of scattering matrices as in \cite{martinssonFastDirectSolvers2019}, 
or methods on graphs \cite{spielman2012nearlylineartimealgorithmspreconditioning,trottenberg2000multigrid,koutis2010approachingoptimalitysolvingsdd}. 
\end{remark}

\section{Numerical examples\label{sec:numericalresults}}

\begin{figure}[ht]
  \centering
  \includegraphics[width=0.8\textwidth]{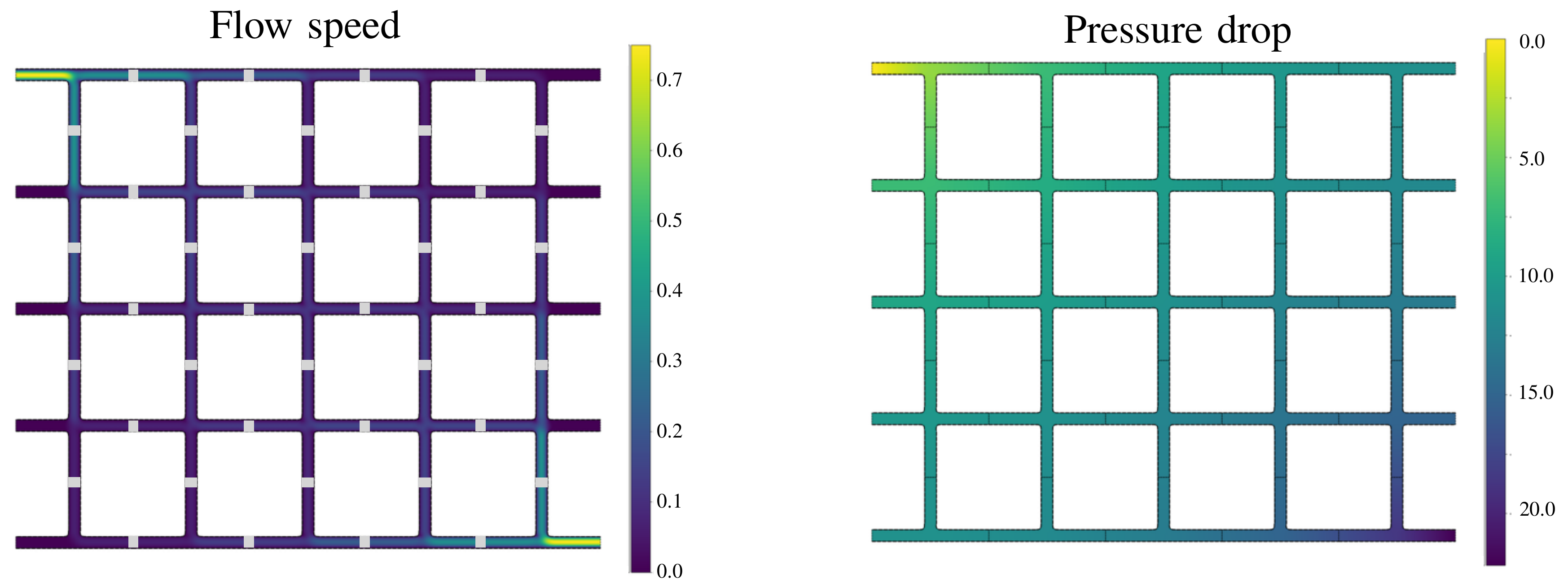}
  \caption{
    Determining the pressure drops needed to establish a desired design flow in a
    multiply-connected geometry. The network is made from 
    25 pieces using 2 distinct component types with 40 junction points, indicated
    with thick, light grey marks on the left.
    The left-hand image shows the flow profile when the inflow was
    set to one at the top left port and the outflow was set to one at the lower right
    port, while the fluxes are set to zero at all other ports.
    In (b), we show the corresponding pressure drop with respect to its value
    at the top left port.
    The scattering matrices for the standard components were computed using an
    iterative solver for the Sherman-Lauricella equation with sixteenth order accuracy 
    and a GMRES tolerance of $10^{-13}$. For validation,
    the global solution was also computed iteratively, with fast multipole acceleration
    and a GMRES tolerance of $10^{-11}$. 
    The estimated relative error in the $L_2$-norm is of the order $10^{-10}$.
    Solving the assembly matrix requires only milliseconds of CPU time.
}
    \label{fig_design}
\end{figure}

To illustrate a typical use case, consider the network illustrated in
\cref{fig_design}, consisting of
25 pieces, made from 2 distinct standard components. 
The goal is to establish a flow, where the input flux is $1$ in the upper left-hand
port and one wishes to guide the flow to the lower right-hand port, with all other 
port fluxes set to zero.
The scattering matrices for the standard components were computed with 12 digits of 
accuracy, and we compare the solution obtained by our merging procedure 
with the iterative solution of the global Sherman-Lauricella
equation applied to the full geometry. Using a GMRES tolerance of $10^{-11}$
for the global integral equation, the two solutions agree to about 
10 digits of accuracy in the $L_2$ norm. Given the scattering matrices,
the required pressure drops are determined by solving the assembly matrix
using only milliseconds of CPU time.

\begin{figure}[ht]
  \centering
  \includegraphics[width=0.8\textwidth]{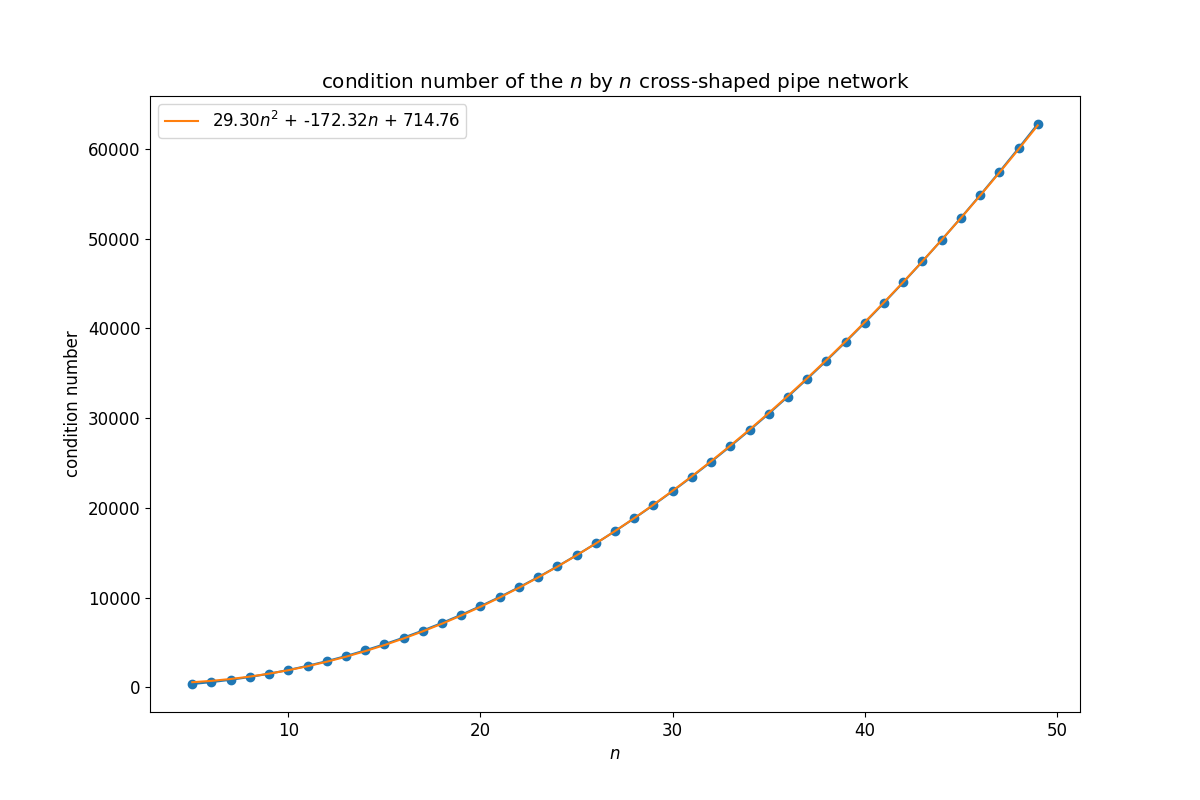}
  \caption{
    The condition number of the assembly matrix for 
    square $n \times n$ networks of the type shown in \cref{fig_design} 
    as a function of $n$.
}
    \label{fig_cond}
\end{figure}

An interesting question concerns the conditioning of the merging process (i.e., 
inversion of the assembly matrix).
While we have not carried out a complete analysis of the problem here, we present the 
results of one experiment.
We consider a network of the type shown in \cref{fig_design}, but with an $n\times n$ grid
replacing the $5\times 5$ grid shown. 
We plot the condition number of the assembly matrix as a function of
$n$ using naive Gaussian elimination in \cref{fig_cond}. The data is extremely well fit by a quadratic
function of $n$. This quadratic scaling of the condition number is similar to that of a grid of electric resistors: given a grid of $n\times n$ resistors in a homogeneous electric circuit grid, the equations given by Kirchhoff's laws have the same quadratic scaling of condition number $\kappa= O(n^2)$ \cite{doyle2000randomwalkselectricnetworks,saadIterativeMethodsSparse2003}.

With $n = 50$ and approximately $2500$ junctions, inverting the assembly matrix
requires 300 milliseconds using MATLAB on a MacBook Air with an M2 processor.
As noted in \cref{sec:merge},
fast algorithms are available that would make inversion of the assembly matrix
even faster -- close to linear scaling in the number of junction points -- but such 
schemes are likely to be of interest only for extremely large systems. 

%\begin{figure}[ht]
%  \centering
%  \includegraphics[width=\textwidth]{example1_fig.pdf}
%  \caption{
%    Assembling the solution of the Stokes equations in a complex branching geometry 
%    consisting of a union of 22 pieces selected involving 7 distinct standard components.
%    The upper plot shows the magnitude of the velocity field in the domain.
%    The black lines mark the locations where the flows are matched.
%    The lower plot shows the absolute error in the velocity field when compared to
%    the solution of the Stokes equations computed in the full geometry.
%    The scattering matrices for the standard components are computed using an
%    iterative solver for the Sherman-Lauricella equation with sixteenth order accuracy 
%    and a GMRES tolerance of $10^{-13}$.
%    The global solution is also computed iteratively, with fast multipole acceleration
%    and a GMRES tolerance of $10^{-11}$. 
%    The estimated relative error in the $L_2$-norm is of the order $10^{-10}$.}
%    \label{fig:connection-error}
%\end{figure}

% Should I include a detailed analysis of the interfacing algorithm here? 
% In short, suppose we can evaluate pressure with accuracy $\epsilon$, 
% then the relative error of fluxes of the generating flows are bounded by m^2\epsilon.
% where m is the total number of generating flows. Better, m^2 can be replaced by the
% number of entries in the matrix generated by interfacing algorithm. 

\section{Conclusions\label{sec:conclusions}}

We have presented a scattering formalism for Stokes flow in branched
networks that reduces the cost of coupling standard components together
to that of solving Kirchhoff's laws in an electrical network.
That is, the coupling involves solving a discrete linear system whose
dimension is the number of locations at which the components are coupled.
The governing PDE is solved only in a pre-computation phase, in order to 
determine the mapping from flow velocities at inlets and outlets to the corresponding
pressure drops, which we refer to as the scattering matrix for the component.
This achieves extremely high accuracy assuming only that -- to either side of the 
coupling -- there is a short length of a straight channel. 
By the nature of slow viscous flow, if the straight channel length is $L$ and the width is $W$, this causes the velocity profile to be 
that of Poiseuille flow at the coupling location with an error that decays 
at the rate  $e^{-4.2 L/W}$ (which we refer to as rapid return to Poiseuille flow).

An important feature of the scattering matrix approach is that the assumption
of separation by many channel widths can be relaxed. To see this, note that
the matching of Poiseuille profiles at the coupling points is simply a way of
enforcing continuity and smoothness of the total flow.
Because of the return to Poiseuille phenomenon, the
velocity is determined by a single constant ($C_0$ in \eqref{papkovich}) as is
the pressure.
Thus, the map from velocity profiles to pressure profiles
has the particularly simple form exploited above.
When components are closer together,
the velocity profile is more complex. Nevertheless, the Sherman-Lauricella equation
can be used to compute a map from the velocity sampled across the inlets and outlets
to a pressure field sampled across the inlets and outlets. 
While the analogy with electrical networks breaks down in this regime,
the essence of the method is unchanged.
Applying the appropriate continuity conditions at the coupling locations,
one constructs a more elaborate assembly matrix.
Extension of our method to this case is under investigation as well as its extension to 
three dimensional flow networks. 

Finally, it is worth noting that
in many applications, the fluid in the network is not homogeneous but rather a colloid,
and accurate flow simulations involve imposing boundary conditions on
the suspended particles. We are currently considering how best to accelerate
such simulations as well.

\section*{Acknowledgements\label{sec:acknowledgements}}

We thank Charlie Peskin, Manas Rachh and Libin Lu for many useful discussions.
The second author gratefully acknowledges support from the Knut and Alice Wallenberg 
foundation under grant 2020.0258.

\end{document}